\documentclass[final]{siamltex}

\usepackage{graphicx}
\usepackage{url}

\title{Stochastic Preconditioning for Iterative Linear Equation Solvers}

\author{Haifeng Qian\thanks{Department of Electrical and Computer Engineering,
University of Minnesota, Minneapolis, Minnesota ({\tt qianhf@ece.umn.edu}).}
\and Sachin S. Sapatnekar\thanks{Department of Electrical and Computer Engineering,
University of Minnesota, Minneapolis, Minnesota ({\tt sachin@ece.umn.edu}).}}

\begin{document}

\maketitle

\begin{abstract}
This paper presents a new stochastic preconditioning approach.
For symmetric diagonally-dominant M-matrices, we prove that
an incomplete LDL factorization can be obtained from random walks,
and used as a preconditioner for an iterative solver,
e.g., conjugate gradient.
It is argued that our factor matrices have better quality,
i.e., better accuracy-size tradeoffs, than preconditioners
produced by existing incomplete factorization methods.
Therefore the resulting preconditioned conjugate gradient (PCG) method
requires less computation than traditional PCG methods
to solve a set of linear equations with the same error tolerance,
and the advantage increases for larger and denser sets of linear equations.
These claims are verified by numerical tests,
and we provide techniques that can potentially extend the theory
to more general types of matrices.
\end{abstract}

\pagestyle{myheadings}
\thispagestyle{plain}
\markboth{H. QIAN AND S. S. SAPATNEKAR}{STOCHASTIC PRECONDITIONING}

\section{Introduction}
\label{chap:intro}

Preconditioning is a crucial part of an iterative linear equation solver.
Suppose a set of linear equations is $A\mathbf{x}=\mathbf{b}$,
where $A$ is a given square nonsingular matrix that is large and sparse,
$\mathbf{b}$ is a given vector,
and $\mathbf{x}$ is the unknown solution vector to be computed.
A (multiplicative) preconditioner is a square nonsingular matrix
$T$ such that an iterative solver can solve the transformed linear
system\footnote{This transformation uses left preconditioning, as against
other options such as
right preconditioning and split preconditioning \cite{saad}.
For simplicity, only left preconditioning is discussed in this paper;
however, the incomplete factorization produced by the proposed approach
can be easily used as a right or split preconditioner.}
$TA\mathbf{x}=T\mathbf{b}$ with a higher convergence rate.

The quality of a preconditioner matrix $T$ is how closely it
approximates\footnote{This can be measured by the spectral radius
of the matrix $\left( I-TA \right)$, or by the condition number
of the matrix $TA$. It is often difficult to analytically quantify
either of the two values for a general matrix, and
the discussion of preconditioner accuracy is mostly qualitative.}
$A^{-1}$.
It is important to note that a preconditioned iterative solver
only requires the evaluation of the product of $T$ and a vector,
and does not require an explicit representation of $T$
in the form of a matrix.
Consequently, any procedure that solves the system $A\mathbf{x}=\mathbf{v}$
approximately can be viewed as a preconditioner, and
the resulting approximate solution can be viewed as the needed product $T\mathbf{v}$,
where $\mathbf{v}$ is any given vector.
Existing preconditioning techniques can be roughly divided into
two categories: explicit methods and implicit methods \cite{benzi}.
In explicit preconditioning methods, which are often referred to
as approximate inverse methods, the preconditioner $T$ is in the form
of a matrix, a product of matrices, or a polynomial of matrices,
and therefore for any given vector $\mathbf{v}$, the product $T\mathbf{v}$
can be evaluated by matrix-vector multiplications \cite{benzi}\cite{saad}.
In implicit preconditioning methods, the preconditioner $T$ is in the form
of $\left( A' \right)^{-1}$, where $A'$ approximates $A$ and is easier to solve,
and therefore for any given vector $\mathbf{v}$, the product $T\mathbf{v}$
is evaluated by solving a linear system with the left-hand-side matrix $A'$ \cite{benzi}\cite{saad}.
Although explicit preconditioning methods have the advantage of being easily parallelizable,
implicit methods have been more successfully developed and more widely used.
A prominent class of implicit preconditioners are those based on incomplete
LU (ILU) factorization; for example, ILU(0), ILU(k) and ILUT are popular
choices in numerical computation\footnote{When matrix $A$ is symmetric and positive definite,
the ILU factors become the corresponding incomplete Cholesky factors, and they
are denoted with the same symbols in this paper.} \cite{templates}\cite{chan}\cite{saad}.

For clarity of the presentation, most of the discussion in this paper is limited to
a special class of left-hand-side matrices: matrix $A$ is referred to as an \emph{R-matrix} if
it is a symmetric M-matrix and is irreducibly diagonally dominant.
One exception is Section~\ref{sec:extension}, which is dedicated to extending
the theory to more general matrix types.

For R-matrices $A$, the most widely used iterative solver is
the Incomplete Cholesky factorization preconditioned Conjugate Gradient (ICCG) method
\cite{templates}\cite{iccg}\cite{saad},
which uses $T = \left( BB^{\rm{T}} \right)^{-1}$ as the preconditioner,
where $B$ is an incomplete Cholesky factor of $A$.
There are various existing techniques to produce $B$ for ICCG. All of these
perform Gaussian elimination on $A$, and use a specific strategy to
drop insignificant entries during the process:
ILU(0) applies a pattern-based strategy, and allows $B_{i,j} \ne 0$ only if $A_{i,j} \ne 0$ \cite{saad};
ILUT applies a value-based strategy, and drops an entry from $B$ if
its value is below a threshold, which is typically determined by 
multiplying the norm of the corresponding row of $A$ by a small constant \cite{saad};
a more advanced strategy can be a combination of pattern, threshold and other
size limits such as maximum number of entries per row.

Our proposed preconditioning technique belongs to the category of multiplicative
implicit preconditioners based on incomplete factorization, and our innovation
is a stochastic procedure for building the incomplete triangular factors.
It is argued theoretically that our factor matrices have better quality,
i.e., better accuracy-size tradeoffs, than preconditioners
produced by existing incomplete factorization methods.
Therefore the resulting preconditioned conjugate gradient (PCG) method,
which we refer to as the \emph{hybrid solver},
requires less computation than traditional PCG methods
to solve a set of linear equations with the same error tolerance,
and the advantage increases for larger and denser sets of linear equations.
We use numerical tests to compare our solver against
both ICCG with ILU(0) and ICCG with ILUT,
and provide techniques that can potentially extend the theory
to more general types of matrices.
Parts of this paper were initially published in \cite{mydac}\cite{hybrid}.

We will now review previous efforts of using stochastic methods
to solve systems of linear equations.
Historically, the theory was developed on two seemingly independent tracks,
related to the analysis of potential theory
\cite{survey1}\cite{antique3}\cite{hersh}\cite{klahr}\cite{survey2}\cite{muller}
and to the solution of systems of linear equations
\cite{antique1}\cite{antique3}\cite{mcsolver2}\cite{mcsolver1}\cite{antique2}.
However, the two applications are closely related and research along each of these
tracks has resulted in the development of analogous algorithms, some of which are
equivalent.

The second of these parallel tracks will be discussed here.
The first work that proposed a random-walk based linear equation solver
is \cite{antique1}, although it was presented as a solitaire game
of drawing balls from urns.
It was proven in \cite{antique1} that, if matrix $A$ satisfies certain conditions,
a game can be constructed and a random variable\footnote{The notations are different from the original ones
used in \cite{antique1}.}
$X$ can be defined such that $E[X]=(A^{-1})_{ij}$, where $(A^{-1})_{ij}$ is
an entry of the inverse matrix of $A$.
Two years later, the work in \cite{antique2} continued this discussion in
the formulation of random walks, and proposed the use of another random 
variable, and it was argued
that, in certain special cases, this variable has a lower variance than $X$, and hence is
likely to converge faster.
Both \cite{antique1} and \cite{antique2} have the advantage of
being able to compute part of an inverse matrix without solving
the whole system, in other words, localizing computation.
Over the years, various descendant stochastic solvers
have been developed \cite{antique3}\cite{mcsolver2}\cite{mcsolver1},
though some of them, e.g., \cite{mcsolver2}\cite{mcsolver1}, do not have the
property of localizing computation.

Early stochastic solvers suffer from accuracy limitations,
and this was remedied by the sequential Monte Carlo method
proposed in \cite{halton62} and \cite{marshall}.
Let $\mathbf{x'}$ be an approximate solution to $A \mathbf{x} = \mathbf{b}$ found by a stochastic solver,
let the residual vector be $\mathbf{r} = \mathbf{b} - A \mathbf{x'}$,
and let the error vector be $\mathbf{z} = \mathbf{x} - \mathbf{x'}$;
then $A \mathbf{z} = \mathbf{r}$.
The idea of the sequential Monte Carlo method is to iteratively solve this
system of equations using a stochastic solver, and in each iteration,
to compute an approximate error vector $\mathbf{z}$ that is
then used to correct the current solution $\mathbf{x'}$.
Although the sequential Monte Carlo method has existed for over forty years,
it has not resulted in any powerful solver that can compete with direct and
iterative solvers, due to the fact that random walks are needed in every
iteration, resulting in a relatively high overall time complexity.

\section{Stochastic Linear Equation Solver}
\label{chap:stoch}

In this section, we study the underlying stochastic mechanism of the
proposed preconditioner.
It is presented as a stand-alone stochastic linear equation solver;
however, in later sections, its usage is not to solve equations,
but to build an incomplete factorization.

\subsection{The Generic Algorithm}
\label{sec:genericgame}

\begin{figure}[htb]
\centering
\includegraphics[width=3in]{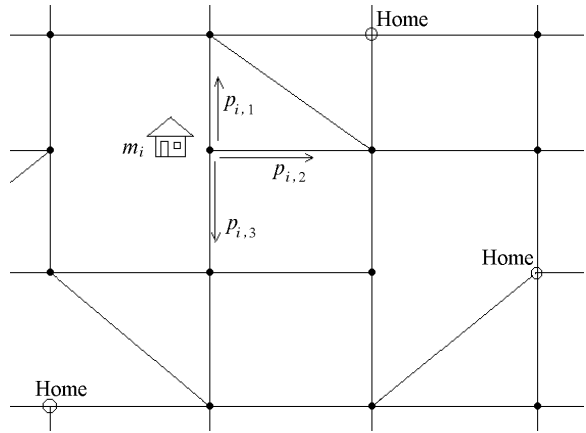}
\caption{An instance of a random walk ``game.''}
\label{dcgame}
\end{figure}

Let us consider a random walk ``game'' defined
on a finite undirected connected graph representing a street map,
for example, Figure~\ref{dcgame}.
A walker starts from one of the nodes,
and every day, he/she goes to an adjacent node $l$ with probability $p_{i,l}$ for $l = 1, 2, \cdots,\mathrm{degree}(i)$,
where $i$ is the current node, $\mathrm{degree}(i)$ is the number of edges connected to node $i$,
and the adjacent nodes are labeled $1, 2, \cdots, \mathrm{degree}(i)$.
The transition probabilities satisfy
\begin{equation}
\sum_{l=1}^{\mathrm{degree}(i)}{p_{i,l}}=1 \label{eq:sumis1}
\end{equation}
The walker pays an amount $m_i$ to a motel for lodging everyday,
until he/she reaches one of the homes, which are a subset of the nodes.
Note that the motel price $m_i$ is a function of his/her current location, node $i$.
The game ends when the walker reaches a home node:
he/she stays there and gets awarded a certain amount of money, $m_0$.
We now consider the problem of calculating the expected amount of money
that the walker has accumulated at the end of the walk,
as a function of the starting node, assuming he/she starts with nothing.
The gain function is therefore defined as
\begin{equation}
f(i)=E[\textrm{total money earned }| \textrm{walk starts at node }i] \label{eq:fdef}
\end{equation}
It is obvious that
\begin{equation}
f(\textrm{one of the homes}) = m_0 \label{eq:fhome}
\end{equation}
For a non-home node $i$, again assuming that the nodes adjacent to $i$ are
labeled 1, 2, $\cdots,\mathrm{degree}(i)$, the $f$ variables satisfy
\begin{equation}
f(i) = \sum_{l=1}^{\mathrm{degree}(i)}{p_{i,l}f(l)}-m_i \label{eq:fequ}
\end{equation}
For a random walk game with $N$ non-home nodes,
there are $N$ linear equations similar to the one above,
and the solution to this set of equations will give the exact values of $f$ at all nodes.

In the above equations obtained from a random walk game,
the set of allowable left-hand-side matrices is a superset of the set
of R-matrices\footnote{A left-hand-side matrix from a random walk game
has all the properties of an R-matrix, except that
it is not necessarily symmetric.}.
In other words, given a set of linear equations $A \mathbf{x} = \mathbf{b}$, where 
$A$ is an R-matrix, we can always construct a random walk game that is 
mathematically equivalent, i.e., such that the $f$ values are the
desired solution $\mathbf{x}$.
To do so, we divide the $i^{\rm th}$ equation 
by $A_{i,i}$ to obtain
\begin{equation}
x_i + \sum_{j \ne i}{\frac{A_{i,j}}{A_{i,i}}x_j} = \frac{b_i}{A_{i,i}}
\end{equation}
\begin{equation}
x_i = \sum_{j \ne i}{\left(-\frac{A_{i,j}}{A_{i,i}}\right) x_j} + \frac{b_i}{A_{i,i}}
\label{eq:zeroaward}
\end{equation}
Equation (\ref{eq:fequ}) and equation (\ref{eq:zeroaward}) have seemingly parallel structures.
Let $N$ be the dimension of matrix $A$, and let us construct a random walk game with $N$ non-home nodes,
which are labeled $1, 2, \cdots, N$.
Due to the properties of an R-matrix, we have
\begin{itemize}
\setlength{\topsep}{0pt}
\setlength{\partopsep}{0pt}
\setlength{\itemsep}{0pt}
\setlength{\parskip}{0pt}
\setlength{\parsep}{0pt}
\item
$\left(-\frac{A_{i,j}}{A_{i,i}}\right)$ is a non-negative value
and can be interpreted as the transition probability of going from node $i$ to node $j$.
\item
$\left( -\frac{b_i}{A_{i,i}} \right)$ can be interpreted as the motel price $m_i$ at node $i$.
\end{itemize}
However, the above mapping is insufficient due to the fact that
condition (\ref{eq:sumis1}) may be broken: the sum of
the $\left(-\frac{A_{i,j}}{A_{i,i}}\right)$ coefficients is not necessarily one.
In fact, because all rows of matrix $A$ are diagonally dominant, 
the sum of the $\left(-\frac{A_{i,j}}{A_{i,i}}\right)$ coefficients is
always less than or equal to one.
Condition (\ref{eq:sumis1}) can be satisfied if we add an extra
transition probability of going from node $i$ to a home node,
by rewriting (\ref{eq:zeroaward}) as
\begin{eqnarray}
& & x_i = \sum_{j \ne i}{\left(-\frac{A_{i,j}}{A_{i,i}}\right) x_j}
+ \frac{\sum_{\forall j}{A_{i,j}}}{A_{i,i}} \cdot m_0 + \frac{b'_i}{A_{i,i}} \nonumber \\
& &\textrm{where  }
b'_i = b_i -  \sum_{\forall j}{A_{i,j}} \cdot m_0 \label{eq:map}
\end{eqnarray}
It is easy to verify that 
$\frac{\sum_{\forall j}{A_{i,j}}}{A_{i,i}}$ is a non-negative value for an R-matrix,
and that the following mapping establishes the equivalence between
equation (\ref{eq:fequ}) and equation (\ref{eq:map}), while satisfying
(\ref{eq:sumis1}) and (\ref{eq:fhome}).
\begin{itemize}
\setlength{\topsep}{0pt}
\setlength{\partopsep}{0pt}
\setlength{\itemsep}{0pt}
\setlength{\parskip}{0pt}
\setlength{\parsep}{0pt}
\item
$\left(-\frac{A_{i,j}}{A_{i,i}}\right)$ is the transition probability of going from node $i$ to node $j$.
\item
$\frac{\sum_{\forall j}{A_{i,j}}}{A_{i,i}}$ is the transition probability
of going from node $i$ to a home node with award $m_0$.
\item
$\left( -\frac{b'_i}{A_{i,i}} \right)$ is the motel price $m_i$ at node $i$.
\end{itemize}
The choice of $m_0$ is arbitrary because $b'_i$ always compensates for
the $m_0$ term in equation (\ref{eq:map}), and in fact $m_0$ can
take different values in (\ref{eq:map}) for different rows $i$.
Therefore the mapping from an equation set to a game is not unique.
A simple scheme can be to let $m_0=0$, and then $m_i=-\frac{b_i}{A_{i,i}}$.

To find $x_i$, the $i^{\rm th}$ entry of solution vector $\mathbf{x}$,
a natural way is to simulate a certain number of random walks from node $i$
and use the average monetary gain
in these walks as the approximated entry value.
If this amount is averaged over a sufficiently large number of walks
by playing the ``game'' a sufficiently large number of times,
then by the Law of Large Numbers \cite{randombook},
an acceptably accurate solution can be obtained.

According to the Central Limit Theorem \cite{randombook},
the estimation error of the above procedure is asymptotically
a zero-mean Gaussian variable with variance inversely proportional to $M$,
where $M$ is the number of walks. Thus there is an accuracy-runtime tradeoff.
In implementation, instead of fixing $M$, one may employ a stopping criterion
driven by a user-specified error margin $\Delta$ and confidence level $\alpha$:
\begin{equation}
P[-\Delta < x'_i - x_i < \Delta] > \alpha \label{eq:delta}
\end{equation}
where $x'_i$ is the estimated $i^{\rm th}$ solution entry from $M$ walks.

\subsection{Two Speedup Techniques}
\label{sec:tricks}

In this section, we propose two new techniques
that dramatically improve the performance of the stochastic solver.
They will play a crucial role in the proposed preconditioning technique.

\subsubsection{Creating Homes}
\label{subsec:newhome}

As discussed in the previous section,
a single entry in the solution vector $\mathbf{x}$ can be evaluated
by running random walks from its corresponding node in the game.
To find the complete solution $\mathbf{x}$,
a straightforward way is to repeat such procedure for every entry.
This, however, is not the most efficient approach,
since much information can be shared between random walks.

We propose a speedup technique by adding the following rule:
after the computation of $x_i$ is finished according to criterion (\ref{eq:delta}),
node $i$ becomes a new home node in the game with an award amount equal to the
estimated value $x'_i$.
In other words, any later random walk that reaches node $i$
terminates, and is rewarded a money amount equal to the assigned $x'_i$.
Without loss of generality, suppose the nodes are processed in the natural ordering
$1, 2, \cdots, N$, then for walks starting from node $k$, the node set $\{1,2,\cdots,k-1\}$
are homes where the walks terminate (in addition to the original homes generated
from the strictly-diagonally-dominant rows of $A$),
while the node set $\{k,k+1,\cdots,N\}$ are motels where the walks pass by.

One way to interpret this technique is by the following observation about
(\ref{eq:fequ}): there is no distinction between the neighboring nodes
that are homes and the neighboring nodes that are motels, and the only
reason that a walk can terminate at a home node is that its $f$ value
is known and is equal to the award.
In fact, any node can be converted to a home node if we know its $f$ value
and assign the award accordingly.
Our new rule is simply utilizing the estimated $x'_i \approx x_i$
in such a conversion.

Another way to interpret this technique is by looking at the source of the value $x'_i$.
Each walk that ends at a new home and obtains such an award is equivalent to an average
of multiple walks, each of which continues walking from there according to
the original game settings.

With this new method, as the computation for the complete solution $\mathbf{x}$ proceeds,
more and more new home nodes are created in the game.
This speeds up the algorithm dramatically,
as walks from later nodes are carried out in a game with a larger and larger
number of homes, and the average number of steps in each walk is reduced.
At the same time, this method helps convergence without increasing $M$,
because, as mentioned earlier, each walk becomes the average of multiple walks.
The only cost\footnote{The cost discussed here is in the context of
the stochastic solver only. For the proposed preconditioner,
this will no longer be an issue.}
is that the game becomes slightly biased when a new home
node is created, due to the fact that the assigned award value is only
an estimate, e.g.~$x'_i \ne x_i$; overall, the benefit of this technique
dominates its cost.

\subsubsection{Bookkeeping}
\label{subsec:bookkeep}

For the same left-hand-side matrix $A$,
traditional direct linear equation solvers are efficient in computing solutions
for multiple right-hand-side vectors after initial matrix factorization,
since only a forward/backward substitution step is required for each additional solve.
Analogous to a direct solver, we propose a speedup mechanism for the stochastic linear
equation solver.

The mechanism is a bookkeeping technique based on the following observation.
In the procedure of constructing a random walk game discussed in Section~\ref{sec:genericgame},
the topology of the game and the transition probabilities are solely determined
by matrix $A$, and hence do not change when the right-hand-side vector
$\mathbf{b}$ changes.
Only motel prices and award values in the game are linked to $\mathbf{b}$.

When solving a set of linear equations with matrix $A$ for the first time,
we create a journey record for every node in the game,
listing the following information.
\begin{itemize}
\setlength{\topsep}{0pt}
\setlength{\partopsep}{0pt}
\setlength{\itemsep}{0pt}
\setlength{\parskip}{0pt}
\setlength{\parsep}{0pt}
\item
For any node $i$, record the number of walks performed from node $i$.
\item
For any node $i$ and any motel node\footnote{The journey record is stored in
a sparse fashion, and a motel $j$ is included only if walks from node $i$ visit $j$ at least once.}
$j$, record the number of times that walks from node $i$ visit node $j$.
\item
For any node $i$ and any home node\footnote{A home node $j$ is included in the journey record
only if at least one walk from $i$ ends at $j$.}
$j$, which can be either an initial home node in the original game or a new home node
created by the technique from Section~\ref{subsec:newhome},
record the number of walks that start from $i$ and end at $j$.
\end{itemize}

Then, if the right-hand-side vector $\mathbf{b}$ changes while the left-hand-side matrix $A$
remains the same, we do not need to perform random walks again.
Instead, we simply use the journey record repeatedly and assume that the walker
takes the same routes, gets awards at the same locations,
pays for the same motels, and only the award amounts and motel prices have been modified.
Thus, after a journey record is created,
new solutions can be computed by some multiplications and additions efficiently.

Practically, this bookkeeping is only feasible after the
technique from Section \ref{subsec:newhome} is in use,
for otherwise the space complexity can be prohibitive for a large matrix.

In the next section, this bookkeeping technique serves as
an important basis of the proposed preconditioner.
There the bookkeeping scheme itself is modified in such a way
that a rigorous proof is presented in Section~\ref{subsec:nonzero}
showing the fact that the space complexity of the modified bookkeeping is upper-bounded
by the space complexity of the matrix factorization
in a direct solver.

\section{Proof of Incomplete Factorization}
\label{sec:proof}

In this section, we build an incomplete LDL factorization of an R-matrix $A$
by extracting information from the journey record of random walks.
The proof is described in two stages: Section \ref{subsec:approxfactor} proves
that the journey record contains an approximate $L$ factor, and then
Section \ref{subsec:nonzero} proves that its non-zero pattern is a subset
of that of the exact $L$ factor.
The formula of the diagonal $D$ factor is derived in Section \ref{subsec:diagonal}.

The procedure of finding this factorization is independent of
the right-hand-side vector $\mathbf{b}$.
Any appearance of $\mathbf{b}$ in this section is symbolic:
its entries do not participate in actual computation, and the involved equations are true for
any possible $\mathbf{b}$.

\subsection{The Approximate Factorization}
\label{subsec:approxfactor}

Suppose the dimension of matrix $A$ is $N$, and its $k^{\rm th}$ row
corresponds to node $k$ in Figure~\ref{dcgame}, $k=1,2,\cdots,N$.
Without loss of generality, assume that in the stochastic solution,
the nodes are processed in the natural ordering $1,2,\cdots,N$.
According to the speedup technique in Section~\ref{subsec:newhome},
for random walks that start from node $k$, the nodes in the set $\{ 1,2,\cdots,k-1 \}$
are already solved and they now serve as home nodes where
a random walk ends.  The awards for reaching nodes $\{ 1,2,\cdots,k-1 \}$ 
are the estimated values of $\{x_1,x_2,\cdots,x_{k-1} \}$ respectively.
Suppose that in equation (\ref{eq:map}), we choose $m_0 = 0$,
and hence the motel prices are given by $m_i=-\frac{b_i}{A_{i,i}}$,
for $i=k,k+1,\cdots,N$. Further,
\begin{itemize}
\setlength{\topsep}{0pt}
\setlength{\partopsep}{0pt}
\setlength{\itemsep}{0pt}
\setlength{\parskip}{0pt}
\setlength{\parsep}{0pt}
\item
Let $M_k$ be the number of walks carried out from node $k$.
\item
Let $H_{k,i}$ be the number of walks that start from node $k$ and end at node 
$i \in \{ 1,2,\cdots,k-1 \}$.
\item
Let $J_{k,i}$ be the number of times that walks 
from node $k$ pass the motel at node $i \in \{ k,k+1,\cdots,N \}$.
\end{itemize}

Taking the average of the results of the $M_k$ walks from node $k$,
we obtain the following equation for the estimated solution entry.
\begin{equation}
x'_k = \frac{\sum_{i=1}^{k-1} {H_{k,i} x'_i} + \sum_{i=k}^{N}{J_{k,i} \frac{b_i}{A_{i,i}}}}{M_k}
\end{equation}
where $x'_i$ is the estimated value of $x_i$ for $i \in \{ 1,2,\cdots,k-1 \}$.
Note that the awards received at the initial home nodes are ignored in
the above equation since $m_0 = 0$.
Moving the $H_{k,i}$ terms to the left side, we obtain
\begin{equation}
-\sum_{i=1}^{k-1}{\frac{H_{k,i}}{M_k} x'_i} + x'_k = \sum_{i=k}^{N}{ \frac{J_{k,i}}{M_k A_{i,i}} b_i}
\end{equation}
By writing the above equation for $k=1,2,\cdots,N$, and assembling the
$N$ equations together into a matrix form, we obtain
\begin{equation}
Y \mathbf{x'} = Z \mathbf{b} \label{eq:yzmatrix}
\end{equation}
where $\mathbf{x'}$ is the approximate solution produced by the stochastic solver;
$Y$ and $Z$ are two square matrices of dimension $N$ such that
\begin{eqnarray}
Y_{k,k}& =& 1, \,\,\;\qquad\quad \forall k \nonumber \\
Y_{k,i}& =& - \frac{H_{k,i}}{M_k}, \,\:\quad \forall k>i \nonumber \\
Y_{k,i}& =& 0, \,\,\;\qquad\quad \forall k < i \nonumber \\
Z_{k,i}& =& \frac{J_{k,i}}{M_k A_{i,i}}, \quad \forall k \le i \nonumber \\
Z_{k,i}& =& 0, \,\,\;\qquad\quad \forall k > i \label{eq:yzdetail}
\end{eqnarray}

These two matrices $Y$ and $Z$ are the journey record built by the bookkeeping
technique in Section~\ref{subsec:bookkeep}. Obviously $Y$ is a lower triangular
matrix with unit diagonal entries, $Z$ is an upper triangular matrix,
and their entries are independent of the right-hand-side vector $\mathbf{b}$.
Once $Y$ and $Z$ are built from random walks, given any $\mathbf{b}$,
one can apply equation (\ref{eq:yzmatrix}) and find $\mathbf{x'}$ efficiently
by a forward substitution.

It is worth pointing out the physical meaning of the entries in matrix $Y$:
the negative of an entry, $\left( - Y_{k,i} \right)$, is asymptotically
equal to the probability that a random walk from node $k$ ends at node $i$,
when $M_k$ goes to infinity.
Another property of matrix $Y$ is that the sum of every row is zero,
except for the first row where only the first entry is non-zero.

From equation (\ref{eq:yzmatrix}), we have
\begin{equation}
Z^{-1} Y \mathbf{x'} = \mathbf{b}
\end{equation}
Since the vector $\mathbf{x'}$ in the above equation is an approximate solution
to the original set of equations $A \mathbf{x} = \mathbf{b}$, it follows
that\footnote{For any vector $\mathbf{b}$, we have
$\left( Z^{-1} Y \right)^{-1} \mathbf{b} = \mathbf{x'} \approx \mathbf{x} = A^{-1} \mathbf{b}$.
Therefore, $A \left( Z^{-1} Y \right)^{-1} \mathbf{b} \approx \mathbf{b}$, and then
$ \left( I - A \left( Z^{-1} Y \right)^{-1} \right) \mathbf{b} \approx 0$.
Since this is true for any vector $\mathbf{b}$, it must be true for eigenvectors of
the matrix $\left( I - A \left( Z^{-1} Y \right)^{-1} \right)$, and it follows
that the eigenvalues of the matrix $\left( I - A \left( Z^{-1} Y \right)^{-1} \right)$
are all close to zero. Thus we claim that $Z^{-1} Y  \approx A$.}
\begin{equation}
Z^{-1} Y  \approx A
\label{eq:ul}
\end{equation}
Because the inverse of an upper triangular matrix, $Z^{-1}$, is also upper 
triangular, equation (\ref{eq:ul}) is in the form of an approximate ``UL factorization''
of $A$.  The following definition and lemma present a simple relation
between UL factorization and the more commonly encountered 
LU factorization.

\begin{definition}
The operator ${\rm rev}(\cdot)$ is defined on square matrices as follows:
given matrix $A$ of dimension $N$, ${\rm rev}(A)$ is also
a square matrix of dimension $N$, such that
\begin{displaymath}
{\rm rev}(A)_{i,j}=A_{N+1-i,N+1-j} , \quad \forall i,j \in \{ 1,2,\cdots,N \}
\end{displaymath}
\label{revdef}
\end{definition}

In simple terms, the operator ${\rm rev}(\cdot)$ merely inverts the row and column
ordering of a matrix. A simple example of applying this operator is as follows:
\begin{displaymath}
{\rm rev} \left( \left[
\begin{array}{ccc}
 1 & 2 & 3 \\
 4 & 5 & 6 \\
 7 & 8 & 9 \\
\end{array}
\right]
\right) = \left[
\begin{array}{ccc}
 9 & 8 & 7 \\
 6 & 5 & 4 \\
 3 & 2 & 1 \\
\end{array}
\right]
\end{displaymath}

\newtheorem{mylemma}{lemma}
\begin{mylemma}
Let $A=LU$ be the LU factorization of a square matrix $A$,
then ${\rm rev}(A) = {\rm rev}(L){\rm rev}(U)$ is true and is the UL factorization of ${\rm rev}(A)$.
\label{ullemma}
\end{mylemma}

Lemma~\ref{ullemma} is self-evident, and the proof is omitted.
It states that the reverse-ordering of the LU factors of $A$
are the UL factors of reverse-ordered $A$.

Applying Lemma~\ref{ullemma} on equation (\ref{eq:ul}), we obtain
\begin{equation}
{\rm rev} (Z^{-1}) {\rm rev}(Y)  \approx {\rm rev}(A)
\label{eq:afterev}
\end{equation}
Since $A$ is an R-matrix and is symmetric, ${\rm rev}(A)$ must be also symmetric,
and we can take the transpose of both sides, and have
\begin{equation}
\left( {\rm rev}(Y) \right)^{\rm T}  \left( {\rm rev} (Z^{-1}) \right)^{\rm T} \approx {\rm rev}(A)
\label{eq:lu}
\end{equation}
The above equation has the form of a Doolittle LU factorization \cite{duffbook}:
matrix $\left( {\rm rev}(Y) \right)^{\rm T}$ is lower triangular
with unit diagonal entries; matrix $\left( {\rm rev} (Z^{-1}) \right)^{\rm T}$
is upper triangular.

\begin{mylemma}
The Doolittle LU factorization of a square matrix is unique.
\label{lemma:doolittle}
\end{mylemma}

The proof of Lemma~\ref{lemma:doolittle} is omitted.
Let the exact Doolittle LU factorization of ${\rm rev}(A)$ be ${\rm rev}(A) = L_{{\rm rev}(A)}U_{{\rm rev}(A)}$,
and its exact LDL factorization be
${\rm rev}(A) = L_{{\rm rev}(A)} D_{{\rm rev}(A)} \left( L_{{\rm rev}(A)} \right)^{\rm T}$.
Since (\ref{eq:lu}) is an approximate Doolittle LU factorization of ${\rm rev}(A)$,
while the exact Doolittle LU factorization is unique, it must be true that:
\begin{eqnarray}
\left( {\rm rev}(Y) \right)^{\rm T} & \approx & L_{{\rm rev}(A)} \label{eq:y}\\
\left( {\rm rev} (Z^{-1}) \right)^{\rm T} & \approx & U_{{\rm rev}(A)}
= D_{{\rm rev}(A)} \left( L_{{\rm rev}(A)} \right)^{\rm T} \label{eq:z}
\end{eqnarray}

The above two equations indicate that from the matrix $Y$ built by random walks,
we can obtain an approximation to factor $L_{{\rm rev}(A)}$,
and that the matrix $Z$ contains redundant information.
Section~\ref{subsec:diagonal} shows how to estimate matrix $D_{{\rm rev}(A)}$
utilizing only the diagonal entries of matrix $Z$, and hence the rest of $Z$
is not needed at all.
According to equation (\ref{eq:yzdetail}), matrix $Y$ is the award register
in the journey record and keeps track of end nodes of random walks,
while matrix $Z$ is the motel-expense register and keeps track of all
intermediate nodes of walks.
Therefore matrix $Z$ is the dominant portion of the journey record,
and by removing all of its off-diagonal entries, the modified journey
record is significantly smaller than that in the original bookkeeping technique
from Section~\ref{subsec:bookkeep}.
In fact, an upper bound on the number of non-zero entries in matrix $Y$ is
proven in the next section.

\subsection{The Incomplete Non-zero Pattern}
\label{subsec:nonzero}

The previous section proves that an approximate factorization of an R-matrix $A$
can be obtained by random walks. However, it does not constitute a proof
of incomplete factorization, because an incomplete factorization implies that
the non-zero pattern of the approximate factor must be a subset of the non-zero
pattern of the exact factor.
Such a proof is the task of this section:
to prove that an entry of $\left( {\rm rev}(Y) \right)^{\rm T}$
can be possibly non-zero only if the corresponding entry of $L_{{\rm rev}(A)}$ is non-zero.

For $i \ne j$, the $(i,j)$ entry of $\left( {\rm rev}(Y) \right)^{\rm T}$ is as follows,
after applying Definition~\ref{revdef} and equation (\ref{eq:yzdetail}).
\begin{equation}
\left( \left( {\rm rev}(Y) \right)^{\rm T} \right)_{i,j} = Y_{N+1-j,N+1-i}
= - \frac{H_{N+1-j,N+1-i}}{M_{N+1-j}}
\label{eq:nonzero}
\end{equation}
This value is non-zero if and only if $j<i$ and $H_{N+1-j,N+1-i} > 0$.
In other words, at least one random walk starts from node $(N+1-j)$ and ends at node $(N+1-i)$.

\begin{figure}[htb]
\centering
\includegraphics[width=2.5in]{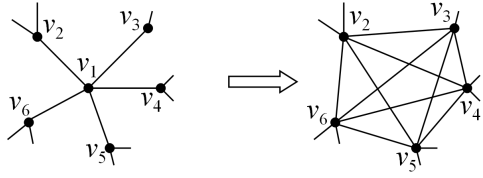}
\caption{One step in Gaussian elimination.}
\label{gauss}
\end{figure}

To analyze the non-zero pattern of $L_{{\rm rev}(A)}$,
certain concepts from the literature of LU factorization are used here,
and certain conclusions are cited without proof.
More details can be found in \cite{amd}\cite{duffbook}\cite{mmd}\cite{directbook}\cite{reachable}.
Figure~\ref{gauss} illustrates one step in the exact Gaussian elimination\footnote{LU factorization
of a matrix is a sequence of Gaussian elimination steps.
From the perspective of the matrix graph, it is a sequence of graph operations that
remove nodes one by one.} of a matrix:
removing one node from the matrix graph, and creating a clique among its neighbors.
For example, when node $v_1$ is removed, a clique is formed for $\{v_2,v_3,v_4,v_5,v_6\}$,
where the new edges correspond to fills added to the remaining matrix.
At the same time, five non-zero values are written into the $L$ matrix,
at the five entries that are the intersections\footnote{In this section,
rows and columns of a matrix are often identified by their corresponding
nodes in the matrix graph, and matrix entries are often identified as
intersections of rows and columns. The reason is that such references are
independent of the matrix ordering, and thereby avoid confusion
due to the two orderings involved in the discussion.}
of node $v_1$'s corresponding column
and the five rows that correspond to nodes $\{v_2,v_3,v_4,v_5,v_6\}$.

\begin{definition}
Given a graph $G=(V,E)$, a node set $S \subset V$, and nodes $v_1,v_2 \in V$ such
that $v_1,v_2 \notin S$,
node $v_2$ is said to be reachable from node $v_1$ through $S$ if there exists a path
between $v_1$ and $v_2$ such that all intermediate nodes, if any, belong to $S$.
\label{reachdef}
\end{definition}

\begin{definition}
Given a graph $G=(V,E)$, a node set $S \subset V$, a node $v_1 \in V$ such that
$v_1 \notin S$, the reachable set of $v_1$ through $S$, denoted 
$R \left( v_1,S \right)$, is defined as:
\begin{displaymath}
R \left( v_1,S \right) = \{v_2 \notin S | v_2 \textrm{ is reachable from } v_1 \textrm{ through } S\}
\end{displaymath}
\end{definition}

Note that if $v_1$ and $v_2$ are adjacent, there is no intermediate node on the path
between them, then Definition~\ref{reachdef} is satisfied, and $v_2$ is reachable from $v_1$
through any node set.
Therefore, $R \left( v_1,S \right)$ always includes the direct neighbors of $v_1$
that do not belong to $S$.

Given an R-matrix $A$, let $G$ be its matrix graph, let $L$ be the complete L factor
in its exact LDL factorization, and let $v_1$ and $v_2$ be two nodes in $G$.
Note that every node in $G$ has a corresponding
row and a corresponding column in $A$ and in $L$.
The following lemma can be derived from \cite{directbook}\cite{reachable}.

\begin{mylemma}
The entry in $L$ at the intersection of column $v_1$ and row $v_2$ is non-zero if and only if:
\begin{enumerate}
\setlength{\topsep}{0pt}
\setlength{\partopsep}{0pt}
\setlength{\itemsep}{0pt}
\setlength{\parskip}{0pt}
\setlength{\parsep}{0pt}
\item
$v_1$ is eliminated prior to $v_2$ during Gaussian elimination
\item
$v_2 \in R \left( v_1, \{ \textrm{nodes eliminated prior to }v_1 \}  \right)$
\end{enumerate}
\end{mylemma}

We now apply this lemma on $L_{{\rm rev}(A)}$.
Because the factorization of ${\rm rev}(A)$ is performed in the reverse ordering, i.e.,
$N,N-1,\cdots,1$, the $(i,j)$ entry of $L_{{\rm rev}(A)}$ is the entry
at the intersection of the column that corresponds to node $(N+1-j)$ and the row
that corresponds to node $(N+1-i)$.
This entry is non-zero if and only if both of the following conditions are met.
\begin{enumerate}
\setlength{\topsep}{0pt}
\setlength{\partopsep}{0pt}
\setlength{\itemsep}{0pt}
\setlength{\parskip}{0pt}
\setlength{\parsep}{0pt}
\item
Node $(N+1-j)$ is eliminated prior to node $(N+1-i)$
\item
$(N+1-i) \in R \left( N+1-j, S_j \right)$\\
where $S_j = \{ \textrm{nodes eliminated prior to }N+1-j \}$
\end{enumerate}

Again, because the Gaussian elimination is carried out in the reverse ordering
$N,N-1,\cdots,1$, the first condition implies that
\begin{eqnarray}
N+1-j & > & N+1-i \nonumber\\
j & < & i \nonumber
\end{eqnarray}
The node set $S_j$ in the second condition is simply
$\{N+2-j,N+3-j,\cdots,N\}$.

Recall that equation (\ref{eq:nonzero}) is non-zero if there is
at least one random walk that starts from node $(N+1-j)$ and ends at 
node $(N+1-i)$. Also recall that according to Section~\ref{subsec:newhome}, when random walks
are performed from node $(N+1-j)$, nodes $\{1,2,\cdots,N-j\}$ are
home nodes that walks terminate, while nodes $S_j = \{N+2-j,N+3-j,\cdots,N\}$ 
are the motel nodes that a walk can pass through.  Therefore, a walk
from node $(N+1-j)$ can possibly end at node $(N+1-i)$,
only if $(N+1-i)$ is reachable from $(N+1-j)$ through the motel node set,
i.e., node set $S_j$.

By now it is proven that both conditions for $\left( L_{{\rm rev}(A)} \right)_{i,j}$ to be
non-zero are necessary conditions for equation (\ref{eq:nonzero}) to be non-zero.
Therefore, the non-zero pattern of $\left( {\rm rev}(Y) \right)^{\rm T}$
is a subset of the non-zero pattern of $L_{{\rm rev}(A)}$.
Together, this conclusion and equation (\ref{eq:y}) give rise to
the following lemma.

\begin{mylemma}
$\left( {\rm rev}(Y) \right)^{\rm T}$ is the $L$ factor of an incomplete LDL factorization of
matrix ${\rm rev}(A)$.
\end{mylemma}

This lemma indicates that, from random walks, we can obtain an incomplete LDL factorization
of the left-hand-side matrix $A$ in its reversed index ordering.
The remaining approximate diagonal matrix $D$ is derived in the next section.

\subsection{The Diagonal Component}
\label{subsec:diagonal}

To evaluate the approximate $D$ matrix, we take the transpose of both
sides of equation (\ref{eq:z}), and obtain
\begin{equation}
{\rm rev} (Z^{-1})  \approx 
 L_{{\rm rev}(A)} D_{{\rm rev}(A)} \label{eq:z2}
\end{equation}

\begin{mylemma}
For a non-singular square matrix $A$, ${\rm rev} (A^{-1}) = \left( {\rm rev}(A) \right) ^{-1}$.
\end{mylemma}

The proof of this lemma is trivial and is omitted.
Applying this lemma on equation (\ref{eq:z2}), we have
\begin{eqnarray}
\left( {\rm rev}(Z) \right) ^{-1} & \approx &
 L_{{\rm rev}(A)} D_{{\rm rev}(A)} \nonumber\\
I & \approx & {\rm rev}(Z) L_{{\rm rev}(A)} D_{{\rm rev}(A)}
\end{eqnarray}
Recall that ${\rm rev}(Z)$ and $L_{{\rm rev}(A)}$ are both lower triangular,
that $L_{{\rm rev}(A)}$ has unit diagonal entries,
and that $D_{{\rm rev}(A)}$ is a diagonal matrix. Therefore, the $( i,i )$ diagonal entry in the
above equation is simply
\begin{eqnarray}
\left( {\rm rev}(Z) \right)_{i,i} \left( L_{{\rm rev}(A)} \right)_{i,i} 
\left( D_{{\rm rev}(A)} \right)_{i,i} & \approx & 1 \nonumber\\
\left( {\rm rev}(Z) \right)_{i,i} \cdot 1 \cdot
\left( D_{{\rm rev}(A)} \right)_{i,i} & \approx & 1 \nonumber\\
\left( D_{{\rm rev}(A)} \right)_{i,i} & \approx & \frac {1} {\left( {\rm rev}(Z) \right)_{i,i}}
\end{eqnarray}
Applying Definition~\ref{revdef} and equation (\ref{eq:yzdetail}), we finally have the equation for computing
the approximate $D$ factor, given as follows.
\begin{eqnarray}
\left( D_{{\rm rev}(A)} \right)_{i,i} & \approx & \frac {1} { Z_{N+1-i,N+1-i}} \nonumber\\
& = & \frac{M_{N+1-i} A_{N+1-i,N+1-i}} {J_{N+1-i,N+1-i}}
\label{eq:d}
\end{eqnarray}

It is worth pointing out the physical meaning of the quantity $\frac {J_{N+1-i,N+1-i}} {M_{N+1-i}} $.
It is the average number of times that a walk from node $N+1-i$ passes node $N+1-i$ itself;
in other words, it is the average number of times that the walker returns to his/her starting point
before the game is over.
Equation (\ref{eq:d}) indicates that an entry in the $D$ factor is equal to the corresponding
diagonal entry of the original matrix $A$ divided by the expected number of returns.

\section{The Hybrid Solver and its Comparison with ILU}
\label{sec:vsilu}

In this section, the proposed hybrid solver for R-matrices is presented in its entirety,
and we argue that it outperforms traditional ICCG methods.

\begin{definition}
The operator ${\rm rev}(\cdot)$ is defined on vectors as follows:
given vector $\mathbf{x}$ of length $N$, ${\rm rev}(\mathbf{x})$ is also
a vector of length $N$, such that
${\rm rev}(\mathbf{x})_i=\mathbf{x}_{N+1-i} ,  \forall i \in \{ 1,2,\cdots,N \}$.
\end{definition}

It is easy to verify that the set of equations $A \mathbf{x} = \mathbf{b}$
is equivalent to
\begin{displaymath}
{\rm rev}(A) {\rm rev}(\mathbf{x}) = {\rm rev}(\mathbf{b})
\end{displaymath}
By now, we have collected the necessary pieces of the proposed hybrid solver, and it
is summarized in the pseudocode in Algorithm~\ref{hybrid}.

\begin{figure}[!ht]
\newtheorem{algorithm}{Algorithm}
\begin{algorithm}
The final hybrid solver for R-matrices:
\em
\small
\begin{center}
\begin{tabular}{l}
\hline
\texttt{Precondition }\{
\\
$\quad$ \texttt{Run random walks, build matrix $Y$ and find diagonal}
\\
$\quad$ $\quad$ $\quad$ \texttt{entries of $Z$ using equation (\ref{eq:yzdetail});}
\\
$\quad$ \texttt{Build $L_{{\rm rev}(A)}$ using equation (\ref{eq:y});}
\\
$\quad$ \texttt{Build $D_{{\rm rev}(A)}$ using equation (\ref{eq:d});}
\\
\}
\\
\texttt{Given $\mathbf{b}$, solve }\{
\\
$\quad$ \texttt{Convert $A \mathbf{x} = \mathbf{b}$ to
${\rm rev}(A) {\rm rev}(\mathbf{x}) = {\rm rev}(\mathbf{b})$;}
\\
$\quad$ \texttt{Apply PCG on ${\rm rev}(A) {\rm rev}(\mathbf{x}) = {\rm rev}(\mathbf{b})$ with the}
\\
$\quad$ $\quad$ $\quad$ \texttt{preconditioner
$\left( L_{{\rm rev}(A)} D_{{\rm rev}(A)} \left( L_{{\rm rev}(A)} \right)^{\rm T} \right) ^{-1}$;}
\\
$\quad$ \texttt{Convert ${\rm rev}(\mathbf{x})$ to $\mathbf{x}$;}
\\
\}
\\
\hline
\end{tabular}
\end{center}
\label{hybrid}
\end{algorithm}
\end{figure}

The proposed hybrid solver essentially replaces the preconditioner in existing ICCG methods
with the incomplete LDL factorization produced by random walks.
We claim that this new preconditioner has better quality than
the incomplete Cholesky factor $B$ produced by traditional incomplete
factorization approaches.
In other words, if matrices $Y$ and $B$ have the same number of non-zero entries,
and given the same target accuracy requirement,
we expect the hybrid solver to converge with fewer iterations than
a traditional ICCG solver preconditioned by $\left( BB^{\rm{T}} \right)^{-1}$.

The argument is based on the fact that, in traditional Gaussian-elimination-based
methods, the operations of eliminating different nodes are correlated and the error
introduced at an earlier node gets propagated to a later node, while in
random walks, the operation on a node is totally independent from other nodes.
We now state this in detail and more precisely.

Let us use the ILUT approach as an example of traditional preconditioning methods;
similar argument can be made for other existing techniques,
as long as they are based on Gaussian elimination.
Suppose in Figure~\ref{gauss}, when eliminating node $v_1$,
the new edge between nodes $v_2$ and $v_3$ corresponds to an entry
whose value falls below a specified threshold,
then ILUT drops that entry from the remaining matrix,
and that edge is removed from the remaining matrix graph.
Later when the algorithm reaches the stage of eliminating node $v_2$, because of that
missing edge, no edge is created from $v_3$ to the neighbors of $v_2$, and thus more edges are
missing, and this new set of missing edges then affect later computations accordingly.
Therefore, an early decision of dropping an entry is propagated throughout the ILUT process.
On the one hand, this leads to the sparsity of $B$, which is desirable; on the other hand,
there is no control over error accumulation, and later columns of $B$ can
deviate from the exact Cholesky factor by an amount that is greater
than the planned threshold of ILUT.
Such error accumulation gets exacerbated for larger and denser matrices.

The hybrid solver does not suffer from this problem.
When we run random walks from node $k$ and collect the $H_{k,i}$ values
to build the $k^{\rm th}$ row of matrix $Y$ according to equation (\ref{eq:yzdetail}),
we only know that the nodes
$\{ 1,2,\cdots,k-1 \}$ are homes, and this is the only information needed.
If, for some reason, the computed $k^{\rm th}$ row of matrix $Y$ is of lower
quality, this error does not affect other rows in any way;
each row is responsible for its own accuracy, according to a criterion to be
discussed in Section~\ref{subsec:stop}.
In fact, in a parallel computing environment, the computation of each row of $Y$
can be assigned to a different processor.

It is worth pointing out that the error accumulation discussed here
is different from the cost of bias discussed at the end of
Section~\ref{subsec:newhome}. That bias in the stochastic solver, in the context
of the hybrid solver, maps to the forward/backward substitution,
i.e., the procedure of applying the preconditioner inside PCG.
Due to the fact that forward/backward substitution is a sequential process,
such bias or error propagation is inevitable in all iterative solvers
as long as an implicit factorization-based multiplicative preconditioner is in use.
Our claim here is that the hybrid solver is free of error accumulation
in building the preconditioner, and not in applying the preconditioner\footnote{After
a row of matrix $Y$ is calculated, it is possible to add a postprocessing
step to drop insignificant entries. The criterion can be any of the strategies
used in traditional incomplete factorization methods, and, as discussed in
Section~\ref{chap:intro}, may be based on pattern, threshold, size limits,
or a combination of them. With such postprocessing, the hybrid solver
still maintains the advantage of independence between row calculations.
This is not included in our implementation.}.

In summary, because of the absence of error accumulation in building the preconditioner,
we expect the hybrid solver to outperform traditional ICCG methods,
and we expect that the advantage becomes more prominent for larger and denser matrices.

\section{Implementation Issues}
\label{sec:implement}

This section describes several implementation aspects of the proposed
preconditioning technique. The goal is twofold:
\begin{itemize}
\setlength{\topsep}{0pt}
\setlength{\partopsep}{0pt}
\setlength{\itemsep}{0pt}
\setlength{\parskip}{0pt}
\setlength{\parsep}{0pt}
\item
To minimize the runtime of building the preconditioner. In other words,
the computation given in the first part of Algorithm~\ref{hybrid}
should be performed with the fewest random walks.
\item
To achieve a better accuracy-size tradeoff.
That is either to improve the accuracy of the preconditioner without increasing
the number of non-zero entries, or to reduce the number of non-zeroes without
losing accuracy\footnote{Again,
the discussion of preconditioner accuracy is mostly qualitative.}.
\end{itemize}

\subsection{Stopping Criterion}
\label{subsec:stop}

The topic of this section is the accuracy control of the preconditioner,
that is, how should one choose $M_k$, the number of walks from node $k$,
to achieve a certain accuracy level
in estimating its corresponding entries in the LDL factorization.
In Section~\ref{sec:genericgame}, the stopping criterion in the stochastic
solver is chosen to be an error margin and a confidence level
defined on the result of a walk;
it is not applicable to the hybrid solver because here it is
necessary for the criterion to be independent of the right-hand-side
vector $\mathbf{b}$.
In our implementation, a new stopping criterion is defined on a value
that is a function of only the left-hand-side matrix $A$, as follows.
Let $\Xi_k = E \left[ \textrm{length of a walk from node }k \right]$,
and let $\Xi'_k$ be the average length of the $M_k$ walks.
The stopping criterion is chosen as
\begin{equation}
P[-\Delta < \frac{\Xi'_k - \Xi_k}{\Xi_k} < \Delta] > \alpha \label{eq:lengthdelta}
\end{equation}
where $\Delta$ is a relative error margin, and $\alpha$ is a confidence level,
for example $\alpha = 99\%$.
Practically, this criterion is checked by the following inequality:
\begin{equation}
\frac{ \Delta \Xi'_k \sqrt{M_k}}{ \sigma_k } > Q^{-1} \left( \frac{1-\alpha}{2} \right)
\end{equation}
where $\sigma_k$ is the standard deviation of the lengths of the $M_k$ walks,
and $Q$ is the standard normal complementary cumulative distribution function.
Thus, $M_k$ is determined dynamically, and random walks are run from node $k$ until
condition (\ref{eq:lengthdelta}) is satisfied.
In practice, it is also necessary to impose a lower bound on $M_k$, e.g., 20 walks.

Note that this is not the only way to design the stopping criterion:
it can also be defined on quantities other than $\Xi_k$ (for example, the
expected number of returns), as long as this quantity does not depend on $\mathbf{b}$.

\subsection{Exact Computations for One-step Walks}
\label{subsec:neighbor}

The implementation technique in this section is a special treatment for the random walks
with length 1, which we refer to as one-step walks.
Such a walk occurs when an immediate neighbor of the starting node is a home node,
and the first step of the walks happens to go there.
The idea is to place stochastic computations performed by one-step walks
with their deterministic limits.

Without loss of generality, assume that the node ordering in the hybrid solver
is the natural ordering $1,2,\cdots,N$.
Let us consider the $M_k$ walks from node $k$, and suppose at least one of its
immediate neighboring nodes
is a home node, which could be either an initial home node
if the $k^{\rm th}$ row of matrix $A$ is strictly diagonally dominant, or a node $j$
such that $j<k$.
Among the $M_k$ walks, let $M_{k,1}$ be the number of one-step walks,
and let $H_{k,i,1}$ be the number of one-step walks that go to node $i$,
where node $i$ is an arbitrary node such that $i<k$.
For the case that node $i$ is not adjacent to node $k$, $H_{k,i,1}$ is simply zero.
For the case that node $i$ is adjacent to node $k$,
note that $H_{k,i,1}$ may not be equal to $M_{k,1}$, as there can be other
immediate neighbors of $k$ that are home nodes.
The $Y_{k,i}$ formula in (\ref{eq:yzdetail}) can be rewritten as
\begin{equation}
Y_{k,i} = - \frac{H_{k,i}}{M_k}
        = - \frac{H_{k,i,1}}{M_k}
          - \left( \frac{M_k - M_{k,1}}{M_k} \right) \cdot \left( \frac{H_{k,i} - H_{k,i,1}}{M_k - M_{k,1}} \right)
\label{eq:yseparate}
\end{equation}

Applying the mapping between transition probabilities and matrix entries in equation (\ref{eq:map}),
the following equations can be derived.
\begin{eqnarray}
\lim_{M_k \to \infty} \frac{H_{k,i,1}}{M_k} & = & P \left[ \textrm{first step goes to node } i \right] \nonumber \\
 & = & - \frac{A_{k,i}}{A_{k,k}}
\end{eqnarray}
\begin{eqnarray}
\lim_{M_k \to \infty} \frac{M_k - M_{k,1}}{M_k} & = & P \left[ \textrm{first step goes to a non-absorbing node} \right] \nonumber \\
 & = & \sum_{j>k}{ P \left[ \textrm{first step goes to node } j \right] } \nonumber \\
 & = & - \frac{ \sum_{j>k}{A_{k,j}} }{A_{k,k}}
\end{eqnarray}

We modify equation (\ref{eq:yseparate}) by replacing the term $\frac{H_{k,i,1}}{M_k}$
and the term $\frac{M_k - M_{k,1}}{M_k}$ with their limits given by the above two equations,
and obtain the following new formula for evaluating $Y_{k,i}$.
\begin{equation}
Y_{k,i} = \frac{A_{k,i}}{A_{k,k}} +
          \left( \frac{ \sum_{j>k}{A_{k,j}} }{A_{k,k}} \right)
          \cdot  \left( \frac{H_{k,i} - H_{k,i,1}}{M_k - M_{k,1}} \right)
\label{eq:newy}
\end{equation}
The remaining stochastic part of this new equation is the term $\frac{H_{k,i} - H_{k,i,1}}{M_k - M_{k,1}}$,
which can be evaluated by considering only random walks whose length is at least two; in other words,
one-step walks are ignored.
In implementation, this can be realized by simulating the first step of walks
by randomly picking one of the non-absorbing neighbors of node $k$; note that
then the number of random walks would automatically be $\left( M_k - M_{k,1} \right)$,
and no adjustment is needed.

With a similar derivation, the $Z_{k,k}$ formula\footnote{Recall that
we only need diagonal entries of matrix $Z$.} in (\ref{eq:yzdetail})
can be modified to
\begin{equation}
Z_{k,k} = \frac{1}{A_{k,k}} + \frac{ \sum_{j>k}A_{k,j} }{A_{k,k}^2}
          - \left( \frac{ \sum_{j>k}{A_{k,j}} }{A_{k,k}^2} \right)
          \cdot \left( \frac{J_{k,k} - J_{k,k,1}}{M_k - M_{k,1}} \right)
\end{equation}
where $J_{k,k,1}$ is the number of times that one-step walks pass node $k$.
Obviously $J_{k,k,1} = M_{k,1}$, and therefore
\begin{equation}
Z_{k,k} = \frac{1}{A_{k,k}} + \left( \frac{ \sum_{j>k}A_{k,j} }{A_{k,k}^2} \right)
          \cdot \left( 1 -  \frac{J_{k,k} - M_{k,1}}{M_k - M_{k,1}} \right)
\label{eq:newz}
\end{equation}
The remaining stochastic part of this new equation, the term $\frac{J_{k,k} - M_{k,1}}{M_k - M_{k,1}}$,
again can be evaluated by considering only random walks with length being at least two.
Practically, such computation is concurrent with evaluating $Y_{k,i}$'s based on
equation (\ref{eq:newy}).

The benefit of replacing (\ref{eq:yzdetail}) with equations
(\ref{eq:newy}) and (\ref{eq:newz}) is twofold:
\begin{itemize}
\setlength{\topsep}{0pt}
\setlength{\partopsep}{0pt}
\setlength{\itemsep}{0pt}
\setlength{\parskip}{0pt}
\setlength{\parsep}{0pt}
\item
Part of the evaluation of $Y_{k,i}$ and $Z_{k,k}$ entries is converted from
stochastic computation to its deterministic limit, and the accuracy is potentially improved.
For the case when all neighbors of node $k$ have lower indices,
i.e., when all neighbors are home nodes,
equations (\ref{eq:newy}) and (\ref{eq:newz}) become exact: they translate to the
exact values of the corresponding entries in the complete LDL factorization.
\item
By avoiding simulating one-step walks, the amount of computation in building the
preconditioner is reduced.
For the case when all neighbors of node $k$ are home nodes, the stochastic parts
of (\ref{eq:newy}) and (\ref{eq:newz}) disappear, and hence no walks are needed.
\end{itemize}

\subsection{Reusing Walks}
\label{subsec:reuse}

Without loss of generality, assume that the node ordering in the hybrid solver
is the natural ordering $1,2,\cdots,N$.
A sampled random walk is completely specified by the node indices along the way,
and hence can be viewed as a sequence of integers
$\{k_1,k_2,\cdots,k_\Gamma\}$,
such that $k_1 > k_\Gamma$, that $k_1 \le k_l, \forall l \in \{ 2,\cdots,\Gamma - 1 \}$,
and that an edge exists between node $k_l$ and node $k_{l+1}$, $\forall l \in \{ 1,\cdots,\Gamma - 1 \}$.
If a sequence of integers satisfy the above requirements,
it is referred to as a legal sequence, and can be mapped to an actual random walk.

Due to the fact that a segment of a legal sequence may also be a legal sequence,
it is possible to extract multiple legal sequences from a single simulated
random walk, and use them also as random walks in the evaluation of
equation (\ref{eq:yzdetail}) or its placement, (\ref{eq:newy}) and (\ref{eq:newz}).
However, there are rules that one must comply with when extracting these legal sequences.
A fundamental premise is that
random samples must be independent of each other.
If two walks share a segment, they become correlated.
Note that if two walks have different starting nodes, they
never participate in the same equation (\ref{eq:newy}) or (\ref{eq:newz}),
and hence are allowed to share segments; if two walks have the same starting nodes,
however, they are prohibited from overlapping.
Moreover, due to the technique in the previous section, any one-step walk should be ignored.

\begin{figure}[htb]
\centering
\begin{displaymath}
\begin{array}{ll}
\textrm{(a)} & \{2,4,6,4,5,7,6,3,2,5,8,1\} \\
\textrm{(b)} & \{4,6,4,5,7,6,3\} \\
             & \{5,7,6,3\} \\
             & \{5,8,1\} \\
\end{array}
\end{displaymath}
\caption[An example of reusing random walks.]
{An example of (a) the legal sequence of a simulated random walk and (b) three extra walks extracted from it.}
\label{sequence}
\end{figure}

Figure~\ref{sequence} shows an example of extracting multiple
legal sequences from a single simulated random walk.
The sequence $\{2,5,8,1\}$ cannot be used because it has the same
starting node as the entire sequence;
the sequence $\{4,5,7,6,3\}$ cannot be used because it has the
same starting node as $\{4,6,4,5,7,6,3\}$ and the two sequences overlap\footnote{It
is also legitimate to extract $\{4,5,7,6,3\}$ instead of $\{4,6,4,5,7,6,3\}$.
However, the premise of random sampling must be fulfilled:
the decision of whether to start a sequence with $k_2=4$ must be made
without the knowledge of numbers after $k_2$,
and the decision of whether to start a sequence with $k_4=4$ must be made
without the knowledge of numbers after $k_4$.
The strategy in Algorithm~\ref{reuse} is to start a sequence as
early as possible, and hence produces $\{4,6,4,5,7,6,3\}$ instead of $\{4,5,7,6,3\}$.}.
On the other hand, $\{5,7,6,3\}$ and $\{5,8,1\}$ are both extracted because
they do not overlap and hence are two independent random walks.

Considering all of the above requirements, the procedure is
shown in Algorithm~\ref{reuse}, where the extracted legal
sequences are directly accounted for in the $M$, $H$ and $J$ accumulators,
which are defined the same as in all equations in this paper.
Note that the simulated random walk is never stored in memory,
and the only extra storage due to this technique is the stacks,
which contain a monotonically increasing sequence of integers at any moment.

\begin{figure}[!ht]
\begin{algorithm}
Extract multiple random walks from a single simulation:
\em
\small
\begin{center}
\begin{tabular}{l}
\hline
\texttt{stack1.push( $k_1$ );}
\\
\texttt{stack2.push( $1$ );}
\\
\texttt{For $l=2,3,\cdots,$ until the end of walk, do }\{
\\
$\quad$ \texttt{While( $k_l<$ stack1.top() )}\{
\\
$\quad$ $\quad$ \texttt{If( $l>$ stack2.top()$+ 1$ )}\{
\\
$\quad$ $\quad$ $\quad$ \texttt{$k'=$ stack1.top();}
\\
$\quad$ $\quad$ $\quad$ \texttt{$M_{k'} = M_{k'} + 1$;}
\\
$\quad$ $\quad$ $\quad$ \texttt{$H_{k',k_l} = H_{k',k_l} + 1$;}
\\
$\quad$ $\quad$ $\quad$ \texttt{$J_{k',k'} = J_{k',k'} + 1$;}
\\
$\quad$ $\quad$ \}
\\
$\quad$ $\quad$ \texttt{stack1.pop();}
\\
$\quad$ $\quad$ \texttt{stack2.pop();}
\\
$\quad$ \}
\\
$\quad$ \texttt{If( $k_l>$ stack1.top() )}\{
\\
$\quad$ $\quad$ \texttt{stack1.push( $k_l$ );}
\\
$\quad$ $\quad$ \texttt{stack2.push( $l$ );}
\\
$\quad$ \}
\\
$\quad$ \texttt{else $J_{k_l,k_l} = J_{k_l,k_l} + 1$;}
\\
\}
\\
\hline
\end{tabular}
\end{center}
\label{reuse}
\end{algorithm}
\end{figure}

This technique reduces the preconditioning runtime by fully utilizing
the information contained in a single simulated random walk,
such that it contributes to equations (\ref{eq:newy}) and (\ref{eq:newz})
as multiple random walks.
It also guarantees that no two overlapping walks have the same starting
node, and hence does not hurt the accuracy of the produced preconditioner.
The only cost of this technique is that the node ordering
of the hybrid solver must be determined beforehand, and hence pivoting
is not allowed during the incomplete factorization\footnote{For R-matrices,
or in general for diagonally dominant matrices, pivoting is not needed.
For more general matrices to be discussed in Section~\ref{sec:extension},
however, the usage of this technique may be limited.}.

\subsection{Matrix Ordering}
\label{subsec:order}

In existing factorization-based preconditioning techniques, matrix ordering
can affect the performance, i.e., the accuracy-size tradeoff, of the preconditioner.
The same statement is true for the proposed stochastic preconditioner.
In general, we perform an incomplete LDL factorization
of the reverse ordering of matrix $A$, we can apply any existing ordering method
on $A$, reverse the ordering that it produces, and then use the resulting ordering.
In this way, any benefit of that ordering method can be inherited by us.
The following are a few examples of practical ordering schemes for
the stochastic preconditioning.
\begin{itemize}
\setlength{\topsep}{0pt}
\setlength{\partopsep}{0pt}
\setlength{\itemsep}{0pt}
\setlength{\parskip}{0pt}
\setlength{\parsep}{0pt}
\item
Approximate minimum degree ordering (AMD) from \cite{amd} is one of the
state-of-the-art ordering techniques to reduce the number of non-zero
entries in a complete LU factorization, or LDL factorization for an R-matrix.
Since the complete L factor has a smaller size, it is likely that with
the same size, the incomplete L factor may have better quality.
Therefore, using a reversed AMD ordering may
improve the accuracy-size tradeoff.
\item
Reverse Cuthill-McKee ordering (RCM) from \cite{cuthillmckee} is a simple
but useful ordering technique to reduce the bandwidth of both the original
matrix $A$ and the complete LU factors, and thereby improve cache efficiency.
The physical CPU time of applying the LU factors on
a particular right-hand-side vector is reduced due to less cache misses.
For the hybrid solver, this means that, with the same preconditioner size,
the actual CPU time of applying the preconditioner may be reduced.
Of course, the ordering to use should be reversed RCM, which becomes
the original Cuthill-McKee ordering.
\item
Random ordering is used in our implementation.
With random ordering, home nodes are relatively evenly distributed at all stages
of the game, and for walks from any node, the most viable home nodes are of
similar distances. Empirically, we have observed a stable performance.
\end{itemize}

\section{Numerical Results}
\label{chap:result}

To evaluate the proposed stochastic preconditioner,
a set of benchmark matrices are generated
by SPARSKIT \cite{sparskit} by finite-difference discretization of the 3D Laplace's
equation $\nabla ^2 u = 0$ with Dirichlet boundary condition. The matrices correspond
to 3D grids with sizes 50-by-50-by-50, 60-by-60-by-60, up to 100-by-100-by-100,
and a right-hand-side vector with all entries being 1 is used with each of them.
They are listed in Table \ref{table1} as benchmarks m1 to m6.
Another four application-specific benchmarks are reported in Table \ref{table2}:
they are placement matrices from VLSI design, and are denser than the
3D-grid matrices.

\begin{table}[ht]
\caption{Computational complexity comparison of the hybrid solver, ICCG with ILU(0) (LASPACK),
and ICCG with ILUT (MATLAB),
to solve for one right-hand-side vector,
for the 3D-grid benchmarks,
with $10^{-6}$ error tolerance.
$N$ is the dimension of a matrix;
$E$ is the number of non-zero entries of a matrix;
$C$ is the number of non-zero entries of the Cholesky factor;
$M1$ is the number of multiplications per iteration;
$I$ is the number of iterations to reach $10^{-6}$ error tolerance;
$M2$ is the total number of multiplications;
$R1$ is the speedup ratio of the hybrid solver over ICCG with ILU(0);
$R2$ is the speedup ratio of the hybrid solver over ICCG with ILUT.}
\label{table1}
\centering
\footnotesize
\tabcolsep 2pt
\begin{tabular}{|l|c|c|c|c|c|c|c|c|c|c|c|c|c|c|c|c|c|}
\hline
Matrix & $N$ & $E$ & \multicolumn{4}{|c|}{ICCG with ILU(0)} & \multicolumn{4}{|c|}{ICCG with ILUT}
& \multicolumn{4}{|c|}{Hybrid} & R1 & R2 \\
\cline{4-17}
 & & & $C$ & $M1$ & $I$ & $M2$ & $C$ & $M1$ & $I$ & $M2$  & $C$ & $M1$ & $I$ & $M2$ & &\\
\hline
m1 & 1.3e5 & 8.6e5 & 4.9e5 & 2.3e6 & 41 & 9.6e7 & 1.7e6 & 4.8e6 & 21 & 1.0e8 & 1.6e6 & 4.5e6 & 18 & 8.1e7 & 1.19 & 1.25 \\
\hline
m2 & 2.2e5 & 1.5e6 & 8.5e5 & 4.1e6 & 48 & 1.9e8 & 3.0e6 & 8.4e6 & 25 & 2.1e8 & 2.8e6 & 7.9e6 & 19 & 1.5e8 & 1.30 & 1.40 \\
\hline
m3 & 3.4e5 & 2.4e6 & 1.4e6 & 6.5e6 & 56 & 3.6e8 & 4.8e6 & 1.3e7 & 29 & 3.9e8 & 4.4e6 & 1.3e7 & 19 & 2.4e8 & 1.51 & 1.61 \\
\hline
m4 & 5.1e5 & 3.5e6 & 2.0e6 & 9.7e6 & 63 & 6.1e8 & 7.2e6 & 2.0e7 & 32 & 6.4e8 & 6.7e6 & 1.9e7 & 19 & 3.6e8 & 1.68 & 1.77 \\
\hline
m5 & 7.3e5 & 5.1e6 & 2.9e6 & 1.4e7 & 71 & 9.8e8 & 1.0e7 & 2.9e7 & 36 & 1.0e9 & 9.6e6 & 2.7e7 & 20 & 5.5e8 & 1.79 & 1.88 \\
\hline
m6 & 1.0e6 & 6.9e6 & 4.0e6 & 1.9e7 & 79 & 1.5e9 & 1.4e7 & 3.9e7 & 40 & 1.6e9 & 1.3e7 & 3.8e7 & 20 & 7.5e8 & 1.99 & 2.09 \\
\hline
\end{tabular}
\end{table}

\begin{table}[ht]
\caption{Computational complexity comparison of the hybrid solver, ICCG with ILU(0) (LASPACK),
and ICCG with ILUT (MATLAB), to solve for one right-hand-side vector,
for the VLSI placement benchmarks, with 1e-6 error tolerance.
$N$, $E$, $C$, $M1$, $I$, $M2$, $R1$ and
$R2$ are as defined in Table \ref{table1}.}
\label{table2}
\centering
\footnotesize
\tabcolsep 2pt
\begin{tabular}{|l|c|c|c|c|c|c|c|c|c|c|c|c|c|c|c|c|c|}
\hline
Matrix & $N$ & $E$ & \multicolumn{4}{|c|}{ICCG with ILU(0)} & \multicolumn{4}{|c|}{ICCG with ILUT}
& \multicolumn{4}{|c|}{Hybrid} & R1 & R2 \\
\cline{4-17}
 & & & $C$ & $M1$ & $I$ & $M2$ & $C$ & $M1$ & $I$ & $M2$  & $C$ & $M1$ & $I$ & $M2$ & &\\
\hline
m7 & 4.3e5 & 5.2e6 & 2.8e6 & 1.3e7 & 122 & 1.5e9 & 6.5e6 & 2.0e7 & 62 & 1.2e9 & 6.5e6 & 2.0e7 & 12 & 2.3e8 & 6.6 & 5.3 \\
\hline
m8 & 3.5e5 & 5.5e6 & 2.9e6 & 1.3e7 & 82 & 1.0e9 & 5.1e6 & 1.7e7 & 27 & 4.6e8 & 5.0e6 & 1.6e7 & 12 & 2.0e8 & 5.3 & 2.4 \\
\hline
m9 & 4.6e5 & 8.2e6 & 4.3e6 & 1.9e7 & 110 & 2.1e9 & 7.5e6 & 2.5e7 & 55 & 1.4e9 & 8.0e6 & 2.5e7 & 13 & 3.3e8 & 6.3 & 4.2 \\
\hline
m10 & 8.8e5 & 9.4e6 & 5.2e6 & 2.3e7 & 159 & 3.7e9 & 1.3e7 & 3.9e7 & 82 & 3.2e9 & 1.2e7 & 3.7e7 & 12 & 4.4e8 & 8.4 & 7.1 \\
\hline
\end{tabular}
\end{table}

\begin{table}[!ht]
\caption {Physical runtimes of the hybrid solver on a Linux workstation with 2.8GHz CPU frequency.
$T1$ is preconditioning CPU time. $T2$ is solving CPU time with 1e-6 error tolerance. Unit is second.}
\label{table4}
\begin{center}
\footnotesize
\begin{tabular}{|l|c|c|c|c|c|c|c|c|c|c|}
\hline
Ckt & m1 & m2 & m3 & m4 & m5 & m6 & m7 & m8 & m9 & m10 \\
\hline
$T1$ & 5.38 & 10.39 & 17.97 & 28.78 & 44.01 & 71.57 & 33.00 & 21.67 & 46.91 & 68.90\\
\hline
$T2$ & 2.52 & 5.80  & 10.59 & 17.50 & 28.61 & 41.07 & 11.90 & 9.73 & 17.07 & 26.09\\
\hline
\end{tabular}
\end{center}
\end{table}

In Table \ref{table1} and Table \ref{table2}, we compare the proposed hybrid solver
against ICCG with ILU(0) and ICCG with ILUT.
The complexity metric is the number of double-precision multiplications needed
at the iterative solving stage for the equation set $A \mathbf{x} = \mathbf{b}$,
in order to converge with an error tolerance of $10^{-6}$.
This error tolerance is defined as:
\begin{equation}
\parallel \mathbf{b} - A \mathbf{x} \parallel_2 \; < \; 10^{-6} \cdot 
\parallel \mathbf{b} \parallel_2
\end{equation}

LASPack \cite{laspack} is used for ICCG with ILU(0),
and MATLAB is used for ICCG with ILUT.
There are three node ordering algorithms available in MATLAB: 
minimum degree ordering (MMD) \cite{mmd}, approximate minimum degree ordering (AMD) \cite{amd},
and reverse Cuthill-McKee ordering (RCM) \cite{cuthillmckee}.
AMD results in the best performance on the benchmarks and is used for all tests.
The dropping threshold of ILUT in MATLAB is tuned, and the accuracy-size tradeoff of the
proposed preconditioner is adjusted, such that the sizes of the Cholesky factors produced
by both methods are similar, i.e., the $C$ values in the tables are close.
For LASPack and MATLAB, the $M1$ values are computed using the following equation.
\begin{equation}
M1 = C \cdot 2 + E + N \cdot 4
\label{eq:m1}
\end{equation}
According to the PCG pseudo codes in \cite{templates} and \cite{saad}, the above equation is the
best possible implementation.
The $M1$ values of the hybrid solver is obtained by a detailed count embedded in its
implementation, and in fact equation (\ref{eq:m1}) is roughly true for the hybrid solver as well.

A clear trend can be observed in Table \ref{table1} and Table \ref{table2} that the larger and denser a matrix is,
the more the hybrid solver outperforms ICCG.
This is consistent with our argument in Section \ref{sec:vsilu}: when
the matrix is larger and denser,
the effect of error accumulation in traditional methods becomes stronger.

The physical runtimes are shown in Table~\ref{table4}.
Admittedly, the preconditioning runtime $T1$
is more than the typical runtime of a traditional incomplete factorization;
however, it is not a large overhead, gets easily amortized over
multiple re-solves, and is worthwhile given the speedup
achieved in the solving stage.

A reference implementation of the stochastic preconditioning for R-matrices,
as well as the hybrid solver, is available to the public \cite{solver}.

\section{Extensions}
\label{sec:extension}

So far the discussion
has been limited to R-matrices. This section presents techniques aimed
at extending the theory to more general matrices, and speculates on
potential challenges in future research on this topic.

\subsection{Asymmetric $A$ Matrices}
\label{subsec:asymmetry}

Let us first remove the symmetry requirement on matrix $A$.
Recall that the construction of the random walk game and the derivation
of equation (\ref{eq:afterev}) does not require $A$ to be symmetric.
Therefore, matrices $Y$ and $Z$ can still be obtained from random walks,
and equation (\ref{eq:afterev}) remains true for an asymmetric matrix $A$.
Suppose ${\rm rev}(A) = L_{{\rm rev}(A)} D_{{\rm rev}(A)} U_{{\rm rev}(A)}$,
where $L_{{\rm rev}(A)}$ is a lower triangular matrix with unit diagonal entries,
$U_{{\rm rev}(A)}$ is an upper triangular matrix with unit diagonal entries,
and $D_{{\rm rev}(A)}$ is a diagonal matrix. This is called the LDU factorization \cite{duffbook},
which is a slight variation of the LU factorization,
and it is easy to show, based on Lemma~\ref{lemma:doolittle}, that the LDU
factorization is also unique.
Substituting the factorization into equation (\ref{eq:afterev}), we have
\begin{equation}
{\rm rev} (Z^{-1}) {\rm rev}(Y) \approx L_{{\rm rev}(A)} D_{{\rm rev}(A)} U_{{\rm rev}(A)}
\end{equation}
Based on the uniqueness of LDU factorization, it must be true that
\begin{eqnarray}
{\rm rev}(Y) & \approx & U_{{\rm rev}(A)} \label{eq:y2}\\
{\rm rev} (Z^{-1}) & \approx & L_{{\rm rev}(A)} D_{{\rm rev}(A)} \label{eq:z3}
\end{eqnarray}
By equation (\ref{eq:y2}), we can approximate $U_{{\rm rev}(A)}$ based on $Y$;
by equation (\ref{eq:z3}), and through the same derivation as in Section~\ref{subsec:diagonal},
we can approximate $D_{{\rm rev}(A)}$ based on the diagonal entries of $Z$.
The remaining question is how to obtain $L_{{\rm rev}(A)}$.

Suppose we construct a random walk game based on $A^{\rm T}$ instead of $A$,
and suppose we obtain matrices $Y_{A^{\rm T}}$ and $Z_{A^{\rm T}}$ based on equation (\ref{eq:yzdetail}).
Then according to equation (\ref{eq:y2}), we have
\begin{equation}
{\rm rev}(Y_{A^{\rm T}}) \approx U_{{\rm rev}(A^{\rm T})}
\label{eq:at}
\end{equation}
where $U_{{\rm rev}(A^{\rm T})}$ is the U factor in the LDU factorization of ${\rm rev}(A^{\rm T})$.
It is easy to derive the following
\begin{equation}
{\rm rev}(A^{\rm T}) = \left( {\rm rev}(A) \right)^{\rm T}
 = \left( U_{{\rm rev}(A)} \right)^{\rm T} D_{{\rm rev}(A)} \left( L_{{\rm rev}(A)} \right)^{\rm T}
\end{equation}
Therefore,
\begin{equation}
L_{{\rm rev}(A^{\rm T})} D_{{\rm rev}(A^{\rm T})} U_{{\rm rev}(A^{\rm T})}
 = \left( U_{{\rm rev}(A)} \right)^{\rm T} D_{{\rm rev}(A)} \left( L_{{\rm rev}(A)} \right)^{\rm T}
\end{equation}
Based on the uniqueness of the LDU factorization, it must be true that
\begin{equation}
\left( L_{{\rm rev}(A)} \right)^{\rm T} = U_{{\rm rev}(A^{\rm T})}
\label{eq:atu}
\end{equation}
By (\ref{eq:at}) and (\ref{eq:atu}), we finally have
\begin{equation}
{\rm rev}(Y_{A^{\rm T}}) \approx \left( L_{{\rm rev}(A)} \right)^{\rm T}
\label{eq:y3}
\end{equation}
In other words, we can approximate $L_{{\rm rev}(A)}$ based on $Y_{A^{\rm T}}$.

In summary, when matrix $A$ is asymmetric, we need to construct two random walk
games for $A$ and $A^{\rm T}$, and then based on the two $Y$ matrices and
the diagonal entries of one of the $Z$ matrices\footnote{Due to the uniqueness
of the LDU factorization, it does not matter the diagonals of which $Z$ are used.},
we can approximate the LDU factorization of ${\rm rev}(A)$ based on equations
(\ref{eq:d}), (\ref{eq:y2}), and (\ref{eq:y3}).
The proof of non-zero pattern is similar to Section~\ref{subsec:nonzero},
and with the same conclusion: the non-zero patterns of the resulting
approximate L and U factors are subsets of those of the exact factors.
Both the time complexity and space complexity of preconditioning become roughly
twice those of the symmetric case: this is the same behavior as a traditional ILU.

\subsection{Random Walk Game with Scaling}
\label{subsec:scaling}

By now, the symmetry restriction on matrix $A$ has been removed,
and the remaining requirements on $A$ are:
the diagonal entries must be positive;
the off-diagonal entries must be negative or zero;
$A$ must irreducibly diagonally dominant, both row-wise and column-wise.

\begin{figure}[!ht]
\centering
\includegraphics[width=3in]{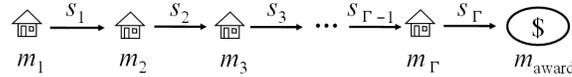}
\caption{A random walk in the modified game with scaling.}
\label{awalk}
\end{figure}

To remove these constraints, a new game is designed by defining a scaling
factor\footnote{A similar concept of scaling factors can be found in
\cite{antique1}, though tailored to its specific game design.}
$s$ on each direction of every edge in the original game from Section~\ref{sec:genericgame}.
Such a scaling factor becomes effective when a random walk passes that particular edge
in that particular direction, and remains effective until this random walk ends.
Let us look at the stochastic solver first.
A walk is shown in Figure~\ref{awalk}: it passes a number of motels, each of
which has its price $m_l$, $l \in \{1,2,\cdots,\Gamma\}$, and ends at a home node with certain award
value $m_{\rm award}$.
The monetary gain of this walk is defined as follows.
\begin{equation}
{\rm gain} = - m_1 - s_1 m_2 - s_1 s_2 m_3 - \cdots
             - \prod_{l=1}^{\Gamma-1}{s_l} \cdot m_\Gamma + \prod_{l=1}^{\Gamma}{s_l} \cdot m_{\rm award}
\label{eq:newgain}
\end{equation}
In simple terms, this new game is different from the original game in that
each transaction amount during the walk is scaled by the product
of the currently active scaling factors.
Define the expected gain function $f$ to be the same as in equation (\ref{eq:fdef}),
and it is easy to derive the replacement of equation (\ref{eq:fequ}):
\begin{equation}
f(i) = \sum_{l=1}^{\mathrm{degree}(i)}{p_{i,l}s_{i,l}f(l)}-m_i
\end{equation}
where $s_{i,l}$ denotes the scaling factor associated with the direction $i \to l$
of the edge between $i$ and $l$, and the rest of the symbols are the same as defined in (\ref{eq:fequ}).

Due to the degrees of freedom introduced by the scaling factors,
the allowable left-hand-side matrix $A$ is now any matrix with non-zero diagonal entries.
In other words, given any matrix $A$ with non-zero diagonal entries, a random
walk game with scaling can be constructed such that the $f$ values, if they
uniquely exist, satisfy a set of linear equations where the left-hand-side matrix is $A$.

A corresponding stochastic preconditioning method can be derived based on this new random walk game,
by redefining the $H$ and $J$ values in equations (\ref{eq:yzdetail}), (\ref{eq:newy}), and (\ref{eq:newz})
to be the sum of products of scaling factors.

If every scaling factor in the game has an absolute value less or equal to 1,
there is no numerical problem in the above new preconditioning procedure.
This can be achieved as long as matrix $A$ is diagonally dominant, in which case
we can simply assign scaling factors to be $+1$ or $-1$, or, if matrix $A$ is complex-valued,
assign complex-valued scaling factors with unit magnitude.
If there exist scaling factors with absolute values over 1, however,
numerical problems may potentially occur since the product of scaling factors may be unbounded.
How to quantify this effect and to analyze the corresponding convergence rate,
is an open question for future research.

Therefore, the conclusion of this section is as follows.
\begin{itemize}
\setlength{\topsep}{0pt}
\setlength{\partopsep}{0pt}
\setlength{\itemsep}{0pt}
\setlength{\parskip}{0pt}
\setlength{\parsep}{0pt}
\item
If the left-hand-side matrix $A$ is irreducibly diagonally dominant both row-wise
and column-wise, the generalized stochastic preconditioning technique
is guaranteed to work, and according to the argument in Section~\ref{sec:vsilu},
the resulting preconditioner
is expected to outperform traditional incomplete factorization methods.
\item
If the left-hand-side matrix $A$ is not diagonally dominant, as long as its
diagonal entries are non-zero, a random walk game exists such that the $f$ values, if they
uniquely exist, satisfy a set of linear equations where the left-hand-side matrix is $A$.
However, no claim is made about the quality of the resulting preconditioner,
and this is open for further investigation.
\end{itemize}

\section*{Acknowledgments}
The authors would like to thank Sani R. Nassif for his contribution
to the stochastic solver, thank Yousef Saad for helpful discussions.


\begin{thebibliography}{99}

\bibitem{amd}
P.~R.~Amestoy, T.~A.~Davis and I.~S.~Duff,
``An approximate minimum degree ordering algorithm,''
\emph{SIAM Journal on Matrix Analysis and Applications},
vol.~17, no.~4, pp.~886-905, 1996.

\bibitem{templates}
R.~Barrett, M.~Berry, T.~F.~Chan, J.~W.~Demmel, J.~Donato, J.~Dongarra,
V.~Eijkhout, R.~Pozo, C.~Romine and H.~A.~van der Vorst,
\emph{Templates for the Solution of Linear Systems: Building Blocks for Iterative Methods},
SIAM, Philadelphia, PA, 1994.

\bibitem{benzi}
M.~Benzi and M.~Tuma,
``A sparse approximate inverse preconditioner for nonsymmetric linear systems,''
\emph{SIAM Journal on Scientific Computing},
vol.~19, no.~3, pp.~968-994, 1998.

\bibitem{chan}
T.~C.~Chan and H.~A.~van der Vorst,
``Approximate and incomplete factorizations,''
Technical Report, Department of Mathematics, University of Utrecht, The Netherlands, 1994.

\bibitem{survey1}
J.~H.~Curtiss,
``Sampling methods applied to differential and difference equations,''
\emph{Proceedings of IBM Seminar on Scientific Computation},
pp.~87-109, 1949.

\bibitem{cuthillmckee}
E.~Cuthill and J.~McKee,
``Reducing the bandwidth of sparse symmetric matrices,''
\emph{Proceedings of the ACM National Conference},
pp.~157-172, 1969.

\bibitem{duffbook}
I.~S.~Duff, A.~M.~Erisman and J.~K.~Reid,
\emph{Direct Methods for Sparse Matrices},
Oxford University Press, New York, NY, 1986.

\bibitem{antique1}
G.~E.~Forsythe and R.~A.~Leibler,
``Matrix inversion by a Monte Carlo method,''
\emph{Mathematical Tables and Other Aids to Computation},
vol.~4, no.~31, pp.~127-129, 1950.

\bibitem{mmd}
A.~George and J.~W.~H.~Liu,
``The evolution of the minimum degree ordering algorithm,''
\emph{SIAM Review},
vol.~31, no.~1, pp.~1-19, 1989.

\bibitem{directbook}
A.~George and J.~W.~H.~Liu,
\emph{Computer Solution of Large Sparse Positive Definite Systems},
Prentice-Hall, Englewood Cliffs, NJ, 1981.

\bibitem{halton62}
J.~H.~Halton,
``Sequential Monte Carlo,''
\emph{Proceedings of the Cambridge Philosophical Society},
vol.~58,~pp.~57-78, 1962.

\bibitem{antique3}
J.~M.~Hammersley and D.~C.~Handscomb,
\emph{Monte Carlo Methods},
Methuen \& Co. Ltd., London, UK, 1964.

\bibitem{reachable}
P.~Heggernes, S.~C.~Eisenstat, G.~Kumfert and A.~Pothen,
``The computational complexity of the Minimum Degree algorithm,''
\emph{Proceedings of 14th Norwegian Computer Science Conference},
pp.~98-109, 2001.

\bibitem{hersh}
R.~Hersh and R.~J.~Griego,
``Brownian motion and potential theory,''
\emph{Scientific American},
vol.~220, pp.~67-74, 1969.

\bibitem{iccg}
D.~S.~Kershaw,
``The incomplete cholesky-conjugate gradient method for the iterative solution of systems of linear equations,''
\emph{Journal of Computational Physics},
vol.~26, pp.~43-65, 1978.

\bibitem{klahr}
C.~N.~Klahr,
``A Monte Carlo method for the solution of elliptic partial differential equations,''
in \emph{Mathematical Methods for Digital Computers}, chap. 14,
John Wiley and Sons, New York, NY, 1962.

\bibitem{survey2}
A.~W.~Knapp,
``Connection between Brownian motion and potential theory,''
\emph{Journal of Mathematical Analysis and Application},
vol.~12, pp.~328-349, 1965.

\bibitem{marshall}
A.~W.~Marshall,
``The use of multi-stage sampling schemes in Monte Carlo,''
\emph{Symposium of Monte Carlo Methods},
pp.~123-140, John Wiley \& Sons, New York, NY, 1956.

\bibitem{muller}
M.~E.~Muller,
``Some continuous Monte Carlo methods for the Dirichlet problem,''
\emph{Annals of Mathematical Statistics},
vol.~27, pp.~569-589, 1956.

\bibitem{mydac}
H.~Qian, S.~R.~Nassif and S.~S.~Sapatnekar,
``Random walks in a supply network,''
\emph{Proceedings of the ACM/IEEE Design Automation Conference},
pp.~93-98, 2003.

\bibitem{hybrid}
H.~Qian and S.~S.~Sapatnekar,
``A hybrid linear equation solver and its application in quadratic placement,''
\emph{ACM/IEEE International Conference on Computer-Aided Design Digest of Technical Papers},
pp.~905-909, 2005.

\bibitem{solver}
H.~Qian and S.~S.~Sapatnekar,
The Hybrid Linear Equation Solver Binary Release. Available at
\url{http://mountains.ece.umn.edu/~sachin/hybridsolver}

\bibitem{saad}
Y.~Saad,
\emph{Iterative Methods for Sparse Linear Systems},
SIAM, Philadelphia, PA, 2003.

\bibitem{sparskit}
Y.~Saad,
SPARSKIT, version 2. Available at\\
\url{http://www-users.cs.umn.edu/~saad/software/SPARSKIT/sparskit.html}

\bibitem{laspack}
T.~Skalicky, LASPack. Available at\\
\url{http://www.mgnet.org/mgnet/Codes/laspack}

\bibitem{mcsolver2}
A.~Srinivasan and V.~Aggarwal,
``Stochastic linear solvers,''
\emph{Proceedings of the SIAM Conference on Applied Linear Algebra}, 2003.

\bibitem{mcsolver1}
C.~J.~K.~Tan and M.~F.~Dixon,
``Antithetic Monte Carlo linear solver,''
\emph{Proceedings of International Conference on Computational Science},
pp.~383-392, 2002.

\bibitem{antique2}
W.~Wasow,
``A note on the inversion of matrices by random walks,''
\emph{Mathematical Tables and Other Aids to Computation},
vol.~6, no.~38, pp.~78-81, 1952.

\bibitem{randombook}
R.~D.~Yates and D.~J.~Goodman,
\emph{Probability and Stochastic Processes: A Friendly Introduction for Electrical and Computer Engineers},
John Wiley \& Sons, New York, NY, 1999.

\end{thebibliography}
\end{document}